\documentclass[12pt]{article}

\usepackage{amsmath,amssymb,a4wide,graphicx}
\usepackage{tikz}
\usepackage[T1]{fontenc}
\usepackage{paralist}
\usepackage{amsthm}
\usepackage{xcolor}
\usepackage{float}
\usepackage[utf8]{inputenc}
\usepackage[english]{babel}
\usepackage{amsmath}
\usepackage{amsfonts}
\usepackage{amssymb}
\usepackage{dsfont}
\usepackage{xcolor}
\usepackage[scale=.75]{geometry}

\newtheorem{theorem}{Theorem}

\newtheorem{coro}{Corollary}

\newcommand{\ml}{l\kern-0.55mm\char39\kern-0.3mm}

\renewcommand{\int}{\mbox{\rm int}}

\def\titlevar{New bounds on domination and independence in graphs}

\usepackage[affil-it]{authblk}

\title{\bf \LARGE \titlevar\footnote{Partially supported by DAAD, Germany (as part of BMBF) and the Ministry of Education, Science, Research and Sport of the Slovak Republic within the project 57447800.}}
\author{Jochen Harant}
\affil{Institut f\"{u}r Mathematik,
TU Ilmenau, 
D-98684 Ilmenau, Germany}
\author{Samuel Mohr\thanks{Gefördert durch die Deutsche Forschungsgemeinschaft (DFG) -- 327533333. This research was supported by the MUNI Award in Science and Humanities of the Grant Agency of Masaryk university.}}
\affil{Faculty of Informatics, Masaryk University, 602 00 Brno, Czech Republic\footnote{Previous affiliation: Institut f\"{u}r Mathematik, TU Ilmenau.}}

\begin{document}

\maketitle

\begin{abstract}
We propose new  bounds  on the domination number and  on the independence number  of a graph and show that our bounds compare favorably to recent ones. Our bounds are obtained by using the {\em Bhatia-Davis} inequality linking the variance, the expected value, the minimum, and the maximum of a random variable with bounded distribution.


\end{abstract}

For a finite set $S$ and an integer $k$,  ${S \choose k}$ and ${|S| \choose k}$ denote the set and the number of $k$-element subsets of $S$, respectively, where, as usual,  ${S \choose k}=\emptyset $ and ${|S| \choose k}=0 $ if $k<0$ or $k>|S|$. Moreover, let $[k]=\{1,2,\cdots,k\}$ if $k\ge 1$.\\
Throughout the paper, we use standard terminology of graph theory and restrict our consideration to the class $\Gamma$ of finite, simple, undirected, connected, and non-complete  graphs. \\\\\\
If we consider a fixed graph $G\in \Gamma$, then let  $V$ be the vertex set of $G$ and $E\subseteq {V \choose 2}$ be the edge set of $G$. Moreover, for $u\in V$, let $N(u)$ be the neighborhood of $u$ in $G$, $d(u)=|N(u)|$ the degree of $u$ in $G$, and $N[u]=N(u)\cup \{u\}$. The minimum degree of $G$ is denoted by  $\delta$. \\
For a positive integer $n$, the class $\Gamma_n$ contains all graphs in $\Gamma$ on $n$ vertices. Clearly, $n\ge 3$ and $\delta\le n-2$ if $G\in \Gamma_n$.\\
A set $D\subseteq V$  is a {\em dominating set} of $G$ if $D\cap N[u]\neq \emptyset$ for every $u\in V$.
The {\em domination number}  $\gamma(G)$ of  $G$ is the minimum cardinality of a dominating set of $G$.\\
We call a set $I\subseteq V$   an {\em independent set} of $G$ if the subgraph $G[I]$ of $G$ induced by $I$ is edgeless.
The {\em independence number}  $\alpha(G)$ of  $G$ is the maximum cardinality of an independent set of $G$.\\
It is hard to determine $\gamma(G)$ or $\alpha(G)$ for $G\in \Gamma$ since the associated decision problems are known to be NP-complete (\cite{GJ}). If we restrict the consideration to bipartite graphs in $\Gamma$, then the computation of $\gamma(G)$ remains hard (\cite{CN}) (note that this is not the case for $\alpha(G)$). For these reasons, we present new upper bounds on the domination number (even for bipartite graphs) in Section~1 and a new lower bound on the independence number in Section~2. All bounds considered here depend on the degrees or even on the edge set of $G\in \Gamma$. In Section~3, our bounds are compared to recent ones. \\
It is well-known how to use the {\em alteration principle} as a powerful tool  to prove upper bounds on $\gamma(G)$ and lower bounds on $\alpha(G)$ (e.g. see \cite{AS}): if $X\subseteq V$ is  randomly chosen, then $D=X\cup Y$ with $Y=\{ u \in V~|~N[u]\cap X=\emptyset\}$ is a dominating set of $G$ and $I=X\setminus Y'$ with $Y'=\{ u \in X~|~N(u)\cap X \neq \emptyset\}$ is an independent set of $G$. Moreover, if $\mu(z)$ and $\mu(z')$ are the expectations of the random variables $z=|X|+|Y|$ and $z'=|X|-|Y'|$, respectively, $min(z)$ is the minimum of $z$, and $max(z')$ denotes the maximum of $z'$, then $\gamma(G)\le min(z)\le \mu(z)$ and $\alpha(G)\ge max(z')\ge \mu(z')$.\\
In \cite{BD}, the {\em Bhatia-Davis} inequality  $Var\le (\mu -min)(max-\mu )$  linking the variance $Var$, the expected value $\mu $, the minimum $min$, and the maximum $max$ of any random
variable with bounded distribution is proved. We use this  result in the proofs of the forthcoming Theorem \ref{dom}, Theorem \ref{dombip}, and Theorem \ref{ind} to strengthen the inequalities $min(z)\le \mu(z)$ and $max(z')\ge \mu(z')$ (see also \cite{ACL}).\\
Our bounds on $\gamma(G)$ and on $\alpha(G)$ are in terms of binomial coefficients. Still, these bounds can be computed efficiently. To see this, we discuss the computational complexity of the factorial. In \cite{HH}, it is shown that the product of two (at most) $s$-digit positive integers can be computed in $O(s\log s)$ time. Combining this  with a result in \cite{SGV}, the computational complexity of  $n!$ is $O(n\log^2 n)$. The slightly weaker statement $O(n(\log n \log \log n)^2)$ can be found in \cite{B}.

\section{Domination}\label{D}

Clark, Shekhtman, Suen, and Fisher \cite{C} proved several upper bounds on $\gamma(G)$ for $G\in \Gamma_n$ in terms of the degrees of $G$, the strongest one being
$\gamma_{CSSF}(G):=\min\limits_{t\in [n-\delta]} \bigg(t+\sum\limits_{u\in V}\frac{{n-d(u)-1 \choose t}}{{n \choose t}}\bigg)$.
We remark that in the definition of $\gamma_{CSSF}(G)$ in \cite{C} the minimum is taken over $t\in [n]$, however, this is not necessary since ${n-d(u)-1 \choose t}=0$ for all $u\in V$ if $t\ge n-\delta$.\\
The forthcoming Theorem \ref{dom} presents an upper bound $\gamma_{HM1}(G)$ on $\gamma(G)$ such that \\
$\gamma_{HM1}(G)\le \gamma_{CSSF}(G)$ for every  $G\in \Gamma$, however, $\gamma_{HM1}(G)$ not only depends on the degrees (as $\gamma_{CSSF}(G)$) but also on the edge set $E$ of $G$. As a consequence of Theorem \ref{dom}, an upper bound $\gamma_{HM2}(G)$ only depending on the degrees of $G$  is presented in Corollary \ref{domcor}.

\begin{theorem}\label{dom}
Let $G\in \Gamma_n$ and, for  $t\in [n-\delta]$,
$a(G,t)=t+\sum\limits_{u \in V}\frac{{n-d(u)-1 \choose t}}{{n \choose t}}$ and \\
$b(G,t)=\sum\limits_{u \in V}\frac{{n-d(u)-1 \choose t}}{{n \choose t}}+2\sum\limits_{\{u,v\} \in { V \choose 2}}\frac{{n-|N[u]\cup N[v]| \choose t}}{{n \choose t}}$.\\
Then $a(G,t)<n$ and $b(G,t)-(a(G,t)-t)^2\ge 0$ for all $t\in [n-\delta]$ and $$\gamma(G)\le \gamma_{HM1}(G):=\min\limits_{ t\in [n-\delta]} \bigg(a(G,t)-\frac{b(G,t)-(a(G,t)-t)^2}{n-a(G,t)}\bigg).$$
\end{theorem}
{\bf Proof.}  Let $ t\in [n-\delta]$ be fixed and pick $X\in {V \choose t}$ uniformly at random.\\
If $Y=\{ u \in V~|~N[u]\cap X=\emptyset\}$, then the set $X\cup Y$ is a dominating set of $G$ and it holds $u \in Y$ if and only if $X\in {V\setminus N[u] \choose t}$.
Consider  $z=|X\cup Y|=|X|+|Y|$ and it follows

$\mu (z)=\mu (|X|)+\mu (|Y|)=t+\sum\limits_{u \in V}P(u \in Y)=t+\sum\limits_{u \in V}P(X\in {V\setminus N[u] \choose t})=t+\sum\limits_{u \in V}\frac{|{V\setminus N[u] \choose t}|}{|{V \choose t}|}$\\
$=t+\sum\limits_{u \in V}\frac{{n-d(u)-1 \choose t}}{{n \choose t}}=a(G,t)$.\\
 Note that $G$ is connected, $|X|<n$, and that $G$ does not contain an edge between $X$ and $Y$. It follows $V\setminus (X\cup Y)\neq \emptyset$,  $z\le n-1$,
 and $a(G,t)=\mu(z)\le max(z)<n$, where $max(z)$ is the maximum of $z$.

If  $min(z)$  denotes   the minimum of $z$,  then $\gamma(G) \le min(z)$ and, by the {\em Bhatia-Davis} inequality,
$\gamma(G) \le \mu (z)-\frac{Var(z)}{n-\mu (z)}=a(G,t)-\frac{Var(z)}{n-a(G,t)}$.\\
Because $Var(z)\ge 0$, it suffices  to show that $Var(z)=b(G,t)-(a(G,t)-t)^2$ to complete the proof.\\
Note that $Var(z)=Var(|X|+|Y|)=Var(t+|Y|)=Var(|Y|)$.\\
If, for $u \in V$,  $Y(u)=1$ if $u \in Y$ and $Y(u)=0$ otherwise, then\\
$Var(|Y|)=Var(\sum\limits_{u \in V}Y(u))=\sum\limits_{u \in V}Var(Y(u))+2\sum\limits_{\{u,v\}\in { V \choose 2}}Cov(Y(u),Y(v))$\\
$=\sum\limits_{u \in V}(\mu (Y(u)^2)-\mu (Y(u))^2)+2\sum\limits_{\{u,v\}\in { V \choose 2}}(\mu (Y(u)Y(v))-\mu (Y(u))\mu (Y(v)))$\\
$=\sum\limits_{u \in V}\mu (Y(u)^2)+2\sum\limits_{\{u,v\}\in { V \choose 2}}\mu (Y(u)Y(v))- (\mu (\sum\limits_{u \in V}Y(u)))^2$\\
$=\sum\limits_{u \in V}\mu (Y(u)^2)+2\sum\limits_{\{u,v\}\in { V \choose 2}}\mu (Y(u)Y(v))- (\mu (|Y|))^2$\\
$=\sum\limits_{u \in V}\mu (Y(u)^2)+2\sum\limits_{\{u,v\}\in { V \choose 2}}\mu (Y(u)Y(v))- (a(G,t)-t)^2$.\\
Note  that $Y(u)Y(v)=1$ if and only if $X\in {V-(N[u]\cup N[v]) \choose t}$. Using $Y(u)^2=Y(u)$, it follows\\
$Var(|Y|)=\sum\limits_{u \in V}\frac{{n-d(u)-1 \choose t}}{{n \choose t}}+2\sum\limits_{\{u,v\}\in { V \choose 2}}\frac{{n-|N[u]\cup N[v]| \choose t}}{{n \choose t}}- (a(G,t)-t)^2=b(G,t)-(a(G,t)-t)^2.$ \hfill $\Box$ \\

Obviously, ${n-|N[u]\cup N[v]| \choose t}\ge {n-d(u)-d(v)-2 \choose t}$ for  $\{u,v\}\in { V \choose 2}$, thus, the forthcoming Corollary~\ref{domcor} is a consequence of Theorem \ref{dom}.

\begin{coro}\label{domcor}
If $G\in \Gamma_n$, $t\in [n-\delta]$,
 $a(G,t)=t+\sum\limits_{u \in V}\frac{{n-d(u)-1 \choose t}}{{n \choose t}}$, and \\$c(G,t)=\sum\limits_{u \in V}\frac{{n-d(u)-1 \choose t}}{{n \choose t}}+2\sum\limits_{\{u,v\}\in { V \choose 2}}\frac{{n-d(u)-d(v)-2 \choose t}}{{n \choose t}}$, then
$$\gamma(G)\le \gamma_{HM2}(G):=\min\limits_{t\in [n-\delta]} (a(G,t)-\frac{c(G,t)-(a(G,t)-t)^2}{n-a(G,t)}).$$
\end{coro}

For a bipartite graph $G\in \Gamma$, several upper bounds on $\gamma(G)$ depending on the cardinalities of the bipartite sets, the size, the minimum degree, or the maximum degree of $G$ can be found in the literature (e.g. see \cite{HP}, \cite{HR}). In Theorem \ref{dombip}, a bound $\gamma_{HM3}(G)$ including much more information on $G$ is presented.

\begin{theorem}\label{dombip}
Let $G$ be a connected bipartite graph with bipartite sets $A$ and $B$ such that $|A|,|B|\ge 2$. \\
For  $(a,b)\in S(|A|,|B|)=\{(a,b)~|~ a\in \{0\}\cup  [|A|],~ b\in \{0\}\cup  [|B|],~ 0<a+b<|A|+|B|\}$, let\\
$e(G,a,b)=a+b+\frac{|A|-a}{|A|}\sum\limits_{u \in A}\frac{{|B|-d(u) \choose b}}{{|B| \choose b}}+\frac{|B|-b}{|B|}\sum\limits_{u \in B}\frac{{|A|-d(u) \choose a}}{{|A| \choose a}}$,\\
$f(G,a,b)=\sum\limits_{u \in A}\bigg(\frac{(|A|- a){|B|-d(u) \choose b}}{|A|{|B| \choose b}}(1-\frac{(|A|- a){|B|-d(u) \choose b}}{|A|{|B| \choose b} })\bigg)$,\\
$g(G,a,b)=\sum\limits_{u \in B}\bigg(\frac{(|B|- b){|A|-d(u) \choose a}}{|B|{|A| \choose a}}(1-\frac{(|B|- b){|A|-d(u) \choose a}}{|B|{|A| \choose a}})\bigg)$,\\
$h(G,a,b)=\sum\limits_{\{u,v\}\in { A \choose 2}}\bigg(\frac{(|A|-a)(|A|-a-1)}{|A|(|A|-1)}\frac{{|B|-|N(u)\cup N(v)| \choose b}}{{|B| \choose b}}-\frac{(|A|- 1)^2}{|A|^2}\frac{{|B|-d(u) \choose b}{|B|-d(v) \choose b}}{{|B| \choose b}^2}\bigg)$, \\
$i(G,a,b)=\sum\limits_{\{u,v\}\in { B \choose 2}}\bigg(\frac{(|B|-b)(|B|-b-1)}{|B|(|B|-1)}\frac{{|A|-|N(u)\cup N(v)| \choose a}}{{|A| \choose a}}-\frac{(|B|- 1)^2}{|B|^2}\frac{{|A|-d(u) \choose a}{|A|-d(v) \choose a}}{{|A| \choose a}^2}\bigg)$, and\\
$j(G,a,b)=\sum\limits_{\{u\}\in A}\sum\limits_{\{v\}\in B}\bigg(\frac{{|A|-|\{u\}\cup N(v)| \choose a}}{{|A| \choose a}}\frac{{|B|-|\{v\}\cup N(u)| \choose b}}{{|B| \choose b}}-\frac{(|A|- a)}{|A|}\frac{{|B|-d(u) \choose b}}{{|B| \choose b}}\frac{(|B|- b)}{|B|}\frac{{|A|-d(v) \choose a}}{{|A| \choose a}}\bigg)$.\\

Then $e(G,a,b)<|A|+|B|$, \\
$k(G,a,b)=f(G,a,b)+g(G,a,b)+2(h(G,a,b)+i(G,a,b)+j(G,a,b))\ge0$, \\
 and $$\gamma(G)\le \gamma_{HM3}(G):=\min\limits_{(a,b)\in S(|A|,|B|)} \bigg(e(G,a,b)-\frac{k(G,a,b)}{|A|+|B|-e(G,a,b)}\bigg)\le \min\{|A|,|B|\}.$$
\end{theorem}
{\bf Proof.}
Let $(a,b)\in S(|A|,|B|)$ be  fixed and pick  $X_A\in {A \choose a}$ and $X_B\in {B \choose b}$ independently and uniformly at random. Clearly, $X_A\cup X_B\neq \emptyset$.\\
If $Y_A=\{ u \in A~|~u\notin X_A, N(u)\cap X_B=\emptyset\}$ and $Y_B=\{ u \in B~|~u\notin X_B, N(u)\cap X_A=\emptyset\}$, then the set $X_A\cup X_B\cup Y_A\cup Y_B$ is a dominating set of $G$. \\
 Clearly, $u \in Y_A$ if and only if $X_A\in {A\setminus \{u\} \choose a}$ and $X_B\in {B\setminus N(u) \choose b}$.
An analog observation holds for $u\in Y_B$.\\
Consider  $z=|X_A\cup X_B\cup Y_A\cup Y_B|=|X_A|+|X_B|+|Y_A|+|Y_B|$ and it follows

$\mu (z)=\mu (|X_A|)+\mu (|X_B|)+\mu (|Y_A|)+\mu (|Y_B|)=a+b+\sum\limits_{u \in A}P(u \in Y_A)+\sum\limits_{u \in B}P(u \in Y_B)$\\
$=a+b+\sum\limits_{u \in A}P(X_A\in {A\setminus \{u\} \choose a})P(X_B\in {B\setminus N(u) \choose b})+\sum\limits_{u \in B}P(X_B\in {B\setminus \{u\} \choose b})P(X_A\in {A\setminus N(u) \choose a})$\\
$=a+b+\frac{{|A|-1 \choose a}}{{|A| \choose a}}\sum\limits_{u \in A}\frac{{|B|-d(u) \choose b}}{{|B| \choose b}}+\frac{{|B|-1 \choose b}}{{|B| \choose b}}\sum\limits_{u \in B}\frac{{|A|-d(u) \choose a}}{{|A| \choose a}}=e(G,a,b)$.\\
If $Y_A\cup Y_B=\emptyset$, then $z=|X_A|+|X_B|=a+b<|A|+|B|$.
If $Y_A\cup Y_B\neq \emptyset$, then $(X_A\cup X_B)\cap(Y_A\cup Y_B)=\emptyset$ and $G$ does not contain an edge between $X_A\cup X_B\neq \emptyset$ and $Y_A\cup Y_B\neq \emptyset$. Since $G$ is connected, we obtain $(A\cup B)\setminus (X_A\cup X_B\cup Y_A\cup Y_B)\neq \emptyset$ and
$z< |A|+|B|$ also in this case. It follows, $e(G,t)=\mu(z)\le max(z)<|A|+|B|$, where $max(z)$ is the maximum of $z$.
If  $min(z)$  is   the minimum of $z$,  then $\gamma(G) \le min(z)$ and
$\gamma(G) \le \mu (z)-\frac{Var(z)}{|A|+|B|-\mu (z)}=e(G,a,b)-\frac{Var(z)}{|A|+|B|-e(G,a,b)}$.\\
We will show that $Var(z)=f(G,a,b)+g(G,a,b)+2(h(G,a,b)+i(G,a,b)+j(G,a,b))$. \\
Because $Var(z)\ge 0$, it follows $k(G,a,b)\ge 0$.\\
Clearly, $Var(z)=Var(|X_A|+|X_B|+|Y_A|+|Y_B|)=Var(a+b+|Y_A|+|Y_B|)$\\
$=Var(|Y_A|+|Y_B|)$.\\
If, for $u \in A$,  $Y_A(u)=1$ if $u \in Y_A$ and $Y_A(u)=0$ otherwise, and, for $u \in B$,  $Y_B(u)=1$ if $u \in Y_B$ and $Y_B(u)=0$ otherwise,  then\\
$Var(z)=Var(\sum\limits_{u \in A}Y_A(u)+\sum\limits_{u \in B}Y_B(u))$\\
$=\sum\limits_{u \in A}Var(Y_A(u))+\sum\limits_{u \in B}Var(Y_B(u))$\\
$+2\bigg(\sum\limits_{\{u,v\}\in { A \choose 2}}Cov(Y_A(u),Y_A(v))+\sum\limits_{\{u,v\}\in { B \choose 2}}Cov(Y_B(u),Y_B(v))+\sum\limits_{\{u\}\in A}\sum\limits_{\{v\}\in B}Cov(Y_A(u),Y_B(v))\bigg)$.\\

We obtain \\
$\sum\limits_{u \in A}Var(Y_A(u))=\sum\limits_{u \in A}(\mu(Y_A(u)^2)-\mu(Y_A(u))^2)=\sum\limits_{u \in A}(\mu(Y_A(u))-\mu(Y_A(u))^2)$\\
$=\sum\limits_{u \in A}P(u\in Y_A)(1-P(u\in Y_A))=\sum\limits_{u \in A}\frac{{|A|-1 \choose a}}{{|A| \choose a}}\frac{{|B|-d(u) \choose b}}{{|B| \choose b}}(1-\frac{{|A|-1 \choose a}}{{|A| \choose a}}\frac{{|B|-d(u) \choose b}}{{|B| \choose b}})=f(G,a,b)$\\
 and, analogously,
$\sum\limits_{u \in B}Var(Y_B(u))=g(G,a,b)$.\\
Moreover,\\
$\sum\limits_{\{u,v\}\in { A \choose 2}}Cov(Y_A(u),Y_A(v))=\sum\limits_{\{u,v\}\in { A \choose 2}}(\mu(Y_A(u)Y_A(v))-\mu(Y_A(u))\mu(Y_A(v)))$\\
$=\sum\limits_{\{u,v\}\in { A \choose 2}}(P(u,v\in Y_A)-P(u\in Y_A)P(v\in Y_A))$\\
$=\sum\limits_{\{u,v\}\in { A \choose 2}}(\frac{{|A|-2 \choose a}}{{|A| \choose a}}\frac{{|B|-|N(u)\cup N(v)| \choose b}}{{|B| \choose b}}-\frac{(|A|- 1)^2}{|A|^2}\frac{{|B|-d(u) \choose b}{|B|-d(v) \choose b}}{{|B| \choose b}^2})=h(G,a,b)$
 and\\
$\sum\limits_{\{u,v\}\in { B \choose 2}}Cov(Y_B(u),Y_B(v))=i(G,a,b)$.\\
Eventually,\\
$\sum\limits_{\{u\}\in A}\sum\limits_{\{v\}\in B}Cov(Y_A(u),Y_B(v))$\\
$=\sum\limits_{\{u\}\in A}\sum\limits_{\{v\}\in B}(P(u\in Y_A, v\in Y_B)-P(u\in Y_A)P( v\in Y_B))=j(G,a,b)$.\\
It remains to show that $\gamma_{HM3}(G)\le \min\{|A|,|B|\}$. If $a=|A|$ and $b=0$, then $X_A=A$, $Y_A=\emptyset$, $X_B=\emptyset$, and, since $G$ is connected,  $Y_B=\emptyset$, thus, $z=|A|$ is a constant. It follows $e(G,|A|,0)=\mu(z)=|A|$, $k(G,|A|,0)=Var(z)=0$, and $\gamma_{HM3}(G)\le |A|$. If $a=0$ and $b=|B|$, then, analogously, $\gamma_{HM3}(G)\le |B|$.
 \hfill $\Box$ \\

We remark that the ideas used in the proofs of Theorem \ref{dom} and Theorem \ref{dombip}  can  be applied  to establish upper bounds on domination numbers of other domination concepts.
As an example, we consider total domination in graphs. A set $D\subseteq V$ of $G$ is {\em total dominating} if $N(u)\cap D \neq \emptyset $ for every $u \in V$ (in contrast to $N[u]\cap D \neq \emptyset $ for a  dominating set $D$). Note that a total dominating set of $G$ exists if and only if  $\delta \ge 1$; $V$ itself is a total dominating set in this case. The minimum cardinality of a total dominating set of $G$ is the {\em total domination number} $\gamma_t(G)$ of $G$. \\
For a randomly chosen set $X\subseteq V$, let $Y=\{u \in V~|~N(u)\cap X \neq \emptyset\}$ and $w_u$ be an arbitrary neighbor of $u \in Y$.  If $S=\{w_u~|~u \in Y\}$, then $X\cup S$ is a total dominating set of $G$, $S\cap X= \emptyset$, and $|S|\le |Y|$. It follows $\gamma_t(G)\le \mu(|X\cup S|)=\mu(|X|+|S|)\le \mu(|X|+|Y|)$. Using the random variable $z=|X|+|Y|$, an upper bound on $\gamma_t(G)$ (as in  Theorem \ref{dom} for $\gamma(t)$) can be proven.
Analog arguments for bipartite graphs lead to a bound on $\gamma_t(G)$ similar to that one for $\gamma(G)$ in Theorem \ref{dombip}.

\section{Independence}\label{I}

 Even though some of the lower bounds on $\alpha(G)$ presented here may also be valid if $G$ is complete or disconnected, we restrict our consideration to $G\in \Gamma$.

Caro and Wei \cite{Caro,Wei} proved the classical lower bound $\alpha_{CW}(G)=\sum\limits_{u \in V}\frac{1}{d(u)+1}$ on  $\alpha(G)$.

It is well-known, that  $\alpha_{CW}(G)>1$ if $G\in \Gamma$.
 The forthcoming bounds $\alpha_{S}(G)$ by Selkow \cite{Selkow} (see also \cite{HM}), $\alpha_{ACL}(G)$ by Angel, Campigotto, and Laforest  \cite{ACL}, and $\alpha_{HR}(G)$  due to Harant and Rautenbach \cite{HR} all strengthen $\alpha_{CW}(G)$:

 $$\alpha(G)\ge \alpha_{S}(G):=\alpha_{CW}(G)+\sum\limits_{u \in V}\frac{1}{d(u)+1}\max\{\frac{d(u)}{d(u)+1}-\sum\limits_{v \in N(u)}\frac{1}{d(v)+1},0 \}.$$

In \cite{ACL}, it is proved that
$A(G)\ge 0$ and that
$\alpha(G)\ge \alpha_{ACL}(G):=\alpha_{CW}(G)+\frac{A(G)}{\alpha_{CW}(G)-1}$
if $G\in \Gamma$ and\\
$A(G)=\sum\limits_{u \in V}\frac{d(u)}{(d(u)+1)^2}-\sum\limits_{\{u,v\} \in E}\frac{2}{(d(u)+1)(d(v)+1)}+\sum\limits_{ \{u,v\} \in {v \choose 2}\setminus E}\frac{2|N(u)\cap N(v)|}{(d(u)+1)(d(v)+1)(2+d(u)+d(v)-|N(u)\cap N(v)|)}$.

The following result is proved in \cite{HR}.\\

{\em
If $G\in \Gamma$,
then there exist
a positive integer $k\in \mathbb{N}$
and   $\psi(u)\in \{0\}\cup [d(u)]$ for all $u \in V$, such that
\begin{align}\label{r1}
\alpha(G)\geq k\geq \sum\limits_{u \in V}\frac{1}{d(u)+1-\psi(u)}
\end{align}
and
\begin{align}\label{r2}\sum\limits_{u \in V}\psi(u)\geq 2(k-1).
\end{align}
}

Consider the function
$\phi(G,l)=\min\{ \sum\limits_{u \in V}\frac{1}{d(u)+1-\psi(u)}~|~\psi(u)\in \{0\}\cup [d(u)],  ~\sum\limits_{u \in V}\psi(u)=l\}$ \\
for $l\in \{0\}\cup [|E|]$.
It follows $\phi(G,0)=\alpha_{CW}(G)$ and, by (\ref{r1}) and (\ref{r2}),  $\alpha(G)\ge \alpha_{HR}(G):=k$ for  $G\in \Gamma$,
where $k$ is the smallest integer such that $k\ge \phi(G,2(k-1))$.
 \\
To determine $\alpha_{HR}(G)$, let ${\cal F}(G,0)$ be the family of values $d(u)+1$ for $u \in V$ (note that a family may contain a member more than once). Moreover, for a non-negative integer $l$, let ${\cal F}(G,l+1)=({\cal F}(G,l)\setminus \{max({\cal F}(G,l))\}) \cup \{max({\cal F}(G,l))-1\}$, where $max({\cal F}(G,l))$ is a maximum member of ${\cal F}(G,l)$. \\
For example,  ${\cal F}(P_3,2)=\{1,2,2\}$ for a path $P_3$ on $3$ vertices.\\
By induction on $l$, it can be seen easily that $\phi(G,l)=\sum\limits_{f\in {\cal F}(G,l)}\frac{1}{f}$ for all $l\in \{0\}\cup [|E|]$, thus, $\alpha_{HR}(G)$ can be calculated by the following algorithm:\\
{\em
1. $k=1$\\
2. \underline{while} $k<\phi(G,2(k-1))$ \underline{do} \\
\underline{begin} $k:=k+1$, go to 2. \underline{end}\\
3. $\alpha_{HR}(G):=k$\\

}

Using the ideas of the proof of Theorem \ref{dom},  the forthcoming Theorem \ref{ind} presents a lower bound $\alpha_{HM}(G)$ on $\alpha(G)$.\\
Note that $\alpha_{S}(G)$, $\alpha_{ACL}(G)$, and $\alpha_{HM}(G)$ depend also on the edge set of $G$ in contrast to the fact that for the calculation of $\alpha_{CW}(G)$ and of $\alpha_{HR}(G)$ only the degrees of $G$ are needed.

\begin{theorem}\label{ind}
Let $G\in \Gamma$ and, for $2\le t\in [n-\delta]$,
$a(G,t)=\sum\limits_{u \in V}{n-d(u)-1 \choose t-1}$ and \\
$b(G,t)=2\sum\limits_{\{u,v\}\in { V \choose 2}\setminus E}\bigg({n-d(u)-1 \choose t-2} +{n-d(v)-1 \choose t-2}-{n-|N[u]\cup N[v]| \choose t-2}\bigg)$.   Then
$$\alpha(G)\ge \alpha_{HM}(G):=\max\limits_{2\le t\in [n-\delta] }\bigg(2t-1-\frac{b(G,t)}{a(G,t)}\bigg).$$
\end{theorem}
{\bf Proof.}
Clearly, $n\ge 3$ and $\delta\le n-2$.
Let $2\le t\in [n-\delta]$ be fixed and pick $X\in {V \choose t}$ uniformly at random.
For convenience, let $a'(G,t)=\frac{a(G,t)}{{n \choose t}}$ and $b'(G,t)=\frac{b(G,t)}{{n \choose t}}$. \\
Note that $P(u \in X)=\frac{t}{n}$ and $P(u,v \in X)=\frac{t(t-1)}{n(n-1)}$ for $u,v \in V, ~u\neq v$.\\
If $Y=\{ u \in X~|~N(u)\cap X \neq \emptyset\}$, then the set $X\setminus Y$ is an independent set of $G$ and it holds $u \in Y$ if and only if $u \in X$ and $X\setminus \{u\}\notin {V\setminus N[u] \choose t-1}$, thus,\\
$P(u \in Y)=P(u \in X \wedge X\setminus \{u\}\notin {V\setminus N[u]\choose t-1})=P(u \in X)\cdot P(X\setminus \{u\}\notin {V\setminus N[u]\choose t-1}~|~u \in X)$\\
$=\frac{t}{n}(1- P(X\setminus \{u\}\in {V\setminus N[u]\choose t-1}~|~u \in X))=\frac{t}{n}(1- \frac{{n-d(u)-1 \choose t-1}}{{n-1 \choose t-1}})$.\\

Consider  $z=|X\setminus Y|=|X|-|Y|\ge 0$ and it follows

$\mu(z)=\mu(|X|)-\mu(|Y|)=t-\sum\limits_{u \in V}P(u \in Y)=t-\sum\limits_{u \in V}\frac{t}{n}(1- \frac{{n-d(u)-1 \choose t-1}}{{n-1 \choose t-1}})=\sum\limits_{u \in V}\frac{{n-d(u)-1 \choose t-1}}{{n \choose t}}=a'(G,t)$ and $\mu(|Y|)=t-a'(G,t)\ge 0$.\\
Since $X=Y$ if  $G[X]$ is connected, it follows $min(z)=0$. Since $t\le n-\delta $, it is possible that $X$ contains a vertex $u \in V$ of minimum degree $\delta$ such that $N(u)\cap X=\emptyset$. \\
In this case, $z\ge 1$, hence, $\mu(z)>0$ for $2\le t\in [n-\delta]$
and we obtain \\
$\alpha(G)\ge \mu(z)+\frac{Var(z)}{\mu(z)}=a'(G,t)+\frac{Var(z)}{a'(G,t)}$ because $\alpha(G)\ge max(z)$.
For $u \in V$, let $Y(u)=1$ if $u \in Y$ and $Y(u)=0$ otherwise. Using $Y(u)^2=Y(u)$, it follows\\
$Var(z)=Var(t-|Y|)=Var(|Y|)=Var(\sum\limits_{u \in V}Y(u))$\\
$=\sum\limits_{u \in V}Var(Y(u))+2\sum\limits_{\{u,v\}\in { V \choose 2}}Cov(Y(u),Y(v))$\\
$=\sum\limits_{u \in V}(\mu(Y(u)^2)-\mu(Y(u))^2)+2\sum\limits_{\{u,v\}\in { V \choose 2}}(\mu(Y(u)Y(v))-\mu(Y(u))\mu(Y(v)))$\\
$=\sum\limits_{u \in V}\mu(Y(u))+2\sum\limits_{\{u,v\}\in { V \choose 2}}\mu(Y(u)Y(v))- \mu(\sum\limits_{u \in V}Y(u))^2$\\
$=\mu(|Y|)+2\sum\limits_{\{u,v\}\in { V \choose 2}}\mu(Y(u)Y(v))- \mu(|Y|)^2$.\\
Let $u\neq v$ and note that  $Y(u)Y(v)=1$ if and only if $X\setminus \{u\}\notin {V\setminus N[u] \choose t-1}$, $X\setminus \{v\}\notin {V\setminus N[v] \choose t-1}$, and $u,v \in X$.\\
If $\{u,v\} \in E$, then $Y(u)Y(v)=1$ if and only if $u,v \in X$, thus, \\
$\mu(Y(u)Y(v))=P(u,v \in X)=P(X\setminus \{u,v\}\in {V\setminus  \{u,v\} \choose t-2})
=\frac{{n-2 \choose t-2}}{{n \choose t}}=\frac{t(t-1)}{n(n-1)}$.\\
If $\{u,v\}\in {V\choose 2}\setminus E$ and $u,v \in Y$, then $X\setminus \{u\}\notin {V\setminus N[u] \choose t-1}$ if and only if $X\setminus \{u,v\}\notin {V\setminus N[u] \choose t-2}$.\\
In this case, $Y(u)Y(v)=1$ if and only if $X\setminus \{u,v\}\notin {V \setminus N[u] \choose t-2}$, $X\setminus \{u,v\}\notin {V\setminus N[v] \choose t-2}$, and $u,v \in X$, thus,\\
$\mu(Y(u)Y(v))=P(Y(u)Y(v)=1)=P(X\setminus \{u,v\}\notin {V\setminus N[u] \choose t-2} \wedge X\setminus \{u,v\}\notin {V\setminus N[v] \choose t-2} \wedge u,v \in X)$\\
$=P(u,v \in X)\cdot P\bigg(X\setminus \{u,v\}\notin {V\setminus N[u] \choose t-2} \wedge X\setminus \{u,v\}\notin {V\setminus N[v] \choose t-2} ~|~ u,v \in X\bigg)$\\
$=\frac{t(t-1)}{n(n-1)}\bigg[1-P\bigg(X\setminus \{u,v\}\in {V\setminus N[u] \choose t-2} \vee X\setminus \{u,v\}\in {V\setminus N[v] \choose t-2} ~|~ u,v \in X\bigg)\bigg]$\\
$=\frac{t(t-1)}{n(n-1)}\bigg[1-P\bigg(X\setminus \{u,v\}\in {V\setminus N[u] \choose t-2}~|~ u,v \in X\bigg) -P\bigg( X\setminus \{u,v\}\in {V\setminus N[v] \choose t-2} ~|~ u,v \in X\bigg)$\\
$+P\bigg(X\setminus \{u,v\}\in {V\setminus N[u] \choose t-2} \wedge X\setminus \{u,v\}\in {V\setminus N[v] \choose t-2} ~|~ u,v \in X\bigg)\bigg]$\\
$=\frac{t(t-1)}{n(n-1)}(1-\frac{{n-d(u)-1 \choose t-2}}{{n-2 \choose t-2}} -\frac{{n-d(v)-1 \choose t-2}}{{n-2 \choose t-2}}+\frac{{n-|N[u]\cup N[v]| \choose t-2}}{{n-2 \choose t-2}})
=\frac{t(t-1)}{n(n-1)}-\frac{{n-d(u)-1 \choose t-2}}{{n \choose t}} -\frac{{n-d(v)-1 \choose t-2}}{{n \choose t}}+\frac{{n-|N[u]\cup N[v]| \choose t-2}}{{n \choose t}}$.\\
We obtain\\
$Var(z)=t-a'(G,t)- (t-a'(G,t))^2+2\sum\limits_{\{u,v\}\in { V \choose 2}}\mu(Y(u)Y(v))$\\
$=t-a'(G,t)- (t-a'(G,t))^2+2\frac{|E|t(t-1)}{n(n-1)}+2\sum\limits_{\{u,v\}\in { V \choose 2}\setminus E}(\frac{t(t-1)}{n(n-1)}-\frac{{n-d(u)-1 \choose t-2}}{{n \choose t}} -\frac{{n-d(v)-1 \choose t-2}}{{n \choose t}}+\frac{{n-|N[u]\cup N[v]| \choose t-2}}{{n \choose t}})$\\
$=t-a'(G,t)- (t-a'(G,t))^2+t(t-1)-2\sum\limits_{\{u,v\}\in { V \choose 2}\setminus E}(\frac{{n-d(u)-1 \choose t-2}}{{n \choose t}} +\frac{{n-d(v)-1 \choose t-2}}{{n \choose t}}-\frac{{n-|N[u]\cup N[v]| \choose t-2}}{{n \choose t}})$.\\
Thus, $\alpha(G)\ge a'(G,t)+\frac{t-a'(G,t)- (t-a'(G,t))^2+t(t-1)-b'(G,t)}{a'(G,t)}=2t-1-\frac{b'(G,t)}{a'(G,t)}=2t-1-\frac{b(G,t)}{a(G,t)}$ for all $t$ such that $2\le t\in [n-\delta]$.\hfill $\Box$

\section{Classification of the Results}\label{C}

In this section, we will compare the upper bounds $\gamma_{CSSF}$, $\gamma_{HM1}$, and $\gamma_{HM2}$  on the domination number  (Theorem \ref{thmCompGamma}, Table \ref{tabGamma}, and Table \ref{tabGammaBip}) and the lower bounds $\alpha_{ACL}$, $\alpha_{S}$, $\alpha_{HR}$, and $\alpha_{HM}$ on the independence number (Theorem \ref{thmCompAlpha}, Table \ref{tabAlpha}, and Table \ref{tabAlphaBip}), where the tables  contain the evaluation   of random graphs.\\
The bound $\gamma_{HM3}$ is not included in this comparison since it is valid for bipartite graphs only and, to our best knowledge, a similar upper bound (also depending on the edge set and being valid for bipartite graphs only) on the domination number of a bipartite graph is not known in the literature. Just an example: for the complete bipartite graph $K_{2,1000}$, it follows $\gamma(K_{2,1000})=\gamma_{HM3}(K_{2,1000})=2$ by Theorem \ref{dombip}, however, computations show that $\gamma_{HM1}(K_{2,1000})>60$ and $\gamma_{HM2}(K_{2,1000}),\gamma_{CSSF}(K_{2,1000})>600$.

We say that two bounds on $\gamma(G)$ or on $\alpha(G)$ are {\em incomparable} if there exist two graphs in $\Gamma$ such that each of the two bounds results in a smaller value than the other bound by one of these graphs.\\

\begin{theorem}\label{thmCompGamma}
$~$
\begin{enumerate}
\item $\gamma_{HM1}(G)\leq \gamma_{CSSF}(G)$ and $\gamma_{HM1}(G)\leq \gamma_{HM2}(G)$ for all graphs $G\in\Gamma$.

\item
$\gamma_{HM2}$ and $\gamma_{CSSF}$ are  incomparable.

\item There is a sequence $\{G_m\}\subseteq\Gamma$ such that $\lim\limits_{m\to\infty}\frac{\gamma_{HM1}(G_m)}{\gamma_{CSSF}(G_m)}=0$.
\end{enumerate}
\end{theorem}

{\bf Proof.}\\
The assertion {\em 1.} of Theorem~\ref{thmCompGamma} follows by Theorem~\ref{dom} and Corollary~\ref{domcor}.

The forthcoming
Table \ref{tabGamma} and Table \ref{tabGammaBip} imply that $\gamma_{HM2}$ and $\gamma_{CSSF}$ are  incomparable.

Eventually, we show that
$\lim\limits_{m\to\infty}\frac{\gamma_{HM1}(G_m)}{\gamma_{CSSF}(G_m)}=0$
for the star $G_m=K_{1,m}$  on $n=m+1\ge 3$ vertices. Actually, we prove $\gamma_{CSSF}(K_{1,m})> \frac{3}{4}m$ and $\gamma_{HM1}(K_{1,m})<2\sqrt{m}$ .\\

For $t\in[m]$, it follows
$a(K_{1,m},t)=t+\sum\limits_{u\in V(K_{1,m})}\frac{\binom{n-d(u)-1}{t}}{\binom{n}{t}}=t+m\cdot \frac{\binom{m-1}{t}}{\binom{m+1}{t}}+0=t+\frac{(m+1-t)(m-t)}{m+1}.$

The  function $t+\tfrac{(m+1-t)(m-t)}{m+1}$  in real $t\in [1,m]$ attains its minimum for $t=\frac{m}{2}$ and we obtain
$\gamma_{CSSF}(K_{1,m})=\min\limits_{ t\in [m]}(a(K_{1,m},t))\ge \frac{m}{2}+\tfrac{(m+2)m}{4(m+1)}> \frac{3}{4}m$.\\

According to Theorem \ref{dom}, let\\
$b(K_{1,m},t)=a(K_{1,m},t)-t+2\sum\limits_{\{u,v\}\in \binom{V(K_{1,m})}{2}}\frac{\binom{n-|N[u]\cup N[v]|}{t}}{\binom{n}{t}}=\frac{(m-t)(m-t+1)}{m+1}+2\binom{m}{2}\cdot \frac{\binom{m-2}{t}}{\binom{m+1}{t}}$\\
$=\frac{(m-t)(m-t+1)}{m+1}+\frac{(m-t+1)(m-t)(m-t-1)}{m+1}=(m-t)\cdot (a(K_{1,m},t)-t)$
and\\
$B(K_{1,m},t):=a(K_{1,m},t)-\frac{b(K_{1,m},t)-(a(K_{1,m},t)-t)^2}{m+1-a(K_{1,m},t)}=a(K_{1,m},t)-\frac{(m-t)\cdot (a(K_{1,m},t)-t)-(a(K_{1,m},t)-t)^2}{(m+1-t)-(a(G_m,t)-t)}$\\
$=t+(a(K_{1,m},t)-t)-\frac{(m-t)\cdot (a(K_{1,m},t)-t)-(a(K_{1,m},t)-t)^2}{(m+1-t)-(a(G_m,t)-t)}=t+\frac{a(K_{1,m},t)-t}{(m+1-t)-(a(G_m,t)-t)}=t+\frac{m-t}{t+1}$\\
$<t+\frac{m}{t+1}$.\\

With $t=\lfloor \sqrt m\rfloor$, we conclude
$\gamma_{HM1}(K_{1,m})=\min\limits_{t\in[m]}B(K_{1,m},t)<\lfloor \sqrt m\rfloor+\frac{m}{\lfloor \sqrt m\rfloor+1}<2\sqrt{m}$.\hfill $\Box$\\ \\

\noindent
\begin{theorem}\label{thmCompAlpha}
$~$
\begin{enumerate}
\item  $\alpha_{ACL}$, $\alpha_{S}$, $\alpha_{HR}$, and $\alpha_{HM}$ are pairwise incomparable.

\item There is a sequence $\{G_m\}\subseteq\Gamma$ such that $\lim\limits_{m\to\infty}\frac{\alpha_{HM}(G_m)}{\max\{\alpha_{ACL}(G_m),\alpha_{HR}(G_m),\alpha_{S}(G_m)\}}=\infty$.\label{alphaSequence}
\end{enumerate}
\end{theorem}

{\bf Proof.}\\
The assertion {\em 1.} of Theorem~\ref{thmCompAlpha} follows from  the forthcoming Table 3 and Table 4.\\
For the proof of {\em 2.} let, for an integer $m\ge 2$, $K_{m,m}\in \Gamma$ be the complete bipartite balanced graph  on $n=2m$ vertices being  $m$-regular. \\
Clearly,  $\alpha_{CW}(K_{m,m})=\alpha_{S}(K_{m,m})=\frac{2m}{m+1}<2$.\\
To calculate $\alpha_{ACL}(K_{m,m})$, note that $A(K_{m,m})=2{m \choose 2}\frac{2m}{(m+1)^2(m+2)}=\frac{2m^2(m-1)}{(m+1)^2(m+2)}$, thus, \\
$\frac{A(K_{m,m})}{\alpha_{CW}(K_{m,m})-1}=\frac{2m^2}{(m+1)(m+2)}<2$, hence,
$\alpha_{ACL}(K_{m,m}):=\alpha_{CW}(K_{m,m})+\frac{A(K_{m,m})}{\alpha_{CW}(K_{m,m})-1}<4$.\\
Recall the definition of $\phi(G,l)$ and it follows $\phi(K_{m,m},2(1-1))=\frac{2m}{m+1}>1$ and\\
$\phi(K_{m,m},2(2-1))=\frac{2m-2}{m+1}+\frac{2}{m}<2$, hence, $\alpha_{HR}(K_{m,m})=2$.\\
Next we will show that $\alpha_{HM}(K_{m,m})\ge (3c-\frac{2c}{1-c})m+o(m)=0.10102...\cdot m+o(m)$, where $c=1-\sqrt{\frac{2}{3}}$.\\
Therefore (see Theorem \ref{ind}), let $f(t)=2t-1-\frac{b(K_{m,m},t)}{a(K_{m,m},t)}$ for $2\le t\in [n-\delta]=[m]$. \\
We obtain $a(K_{m,m},t)=2m{m-1 \choose t-1}$ and
$b(K_{m,m},t)=4{ m \choose 2}\bigg(2{m-1 \choose t-2} -{m-2 \choose t-2}\bigg)$.\\
Using ${m-1 \choose t-2}=\frac{t-1}{m-t+1}{m-1 \choose t-1}$ and ${m-2 \choose t-2}=\frac{t-1}{m-1}{m-1 \choose t-1}$, it follows\\
$f(t)=2t-1-(m-1)(t-1)(\frac{2}{m-t+1}-\frac{1}{m-1})=3t-2-\frac{2(m-1)(t-1)}{m-t+1}>3t-2-2m\frac{t-1}{m-t+1}=g(t)$.\\
We choose $t_0=\lfloor cm\rfloor $ and $s_0=cm-1$, hence, $s_0\le t_0\le s_0+1$.\\
It follows $\alpha_{HM}(K_{m,m})\ge f(t_0)\ge g(t_0)\ge 3s_0-2-\frac{2ms_0}{m-s_0}=(3c-\frac{2c}{1-c})m+o(m)$.\hfill $\Box$\\ \\

We considered two kinds of random graphs. First, let $\mathds{G}(n,p)$ be the probability space of all graphs with vertex set $[n]$ and each possible edge appears with probability $p$ independently from all other choices.
We tested 20,000 connected and non-complete graphs in $\mathds{G}(n,p)$ with 500 graphs in each combination of $p\in\{.2,.3,.5,.6,.8\}$ and $n\in\{10,20,30,50,80,100,120,\\150\}$.\\

Since it seems that our bounds perform well on graphs that have a structure ``close'' to bipartite graphs, we further investigated another probability space $\mathds{B}(n,p_R,p_A)$ of all graphs on $n$ vertices, where $0\leq p_R,p_A\leq 1$.
A graph from $\mathds{B}(n,p_R,p_A)$ can be sampled by starting 
with a complete bipartite graph $K_{k,l}$ on $n=k+l$ vertices for
random $k$ chosen uniformly from
$\{\frac{1}{n},\frac{2}{n}, \dots, \frac{n-1}{n}\}$.
Then each of the $k\cdot l$ edges of $K_{k,l}$ is removed independently
with probability $p_R$ and each of the $\binom k2+\binom l2$
non-edges of $K_{k,l}$ is added independently as an edge
with probability $p_A$.
Note that
$\mathds{B}(n,1-p,p)=\mathds{G}(n,p)$ for all $p$ with $0\le p\le 1$.
We considered 18,000 connected and non-complete graphs with 500 graphs randomly sampled with parameters for each combination of $p_A,p_R\in\{.02,.05,.1\}$ and $n\in\{10,25,50,100\}$.\\

Table~\ref{tabGamma} and Table~\ref{tabGammaBip} should be read as follows: The percentage of considered random graphs describes the cases where the floored integer of the bound on the domination number describing the row is strictly smaller than the floored integer of the bound describing the column.\\

\begin{table}[!ht]\centering
\begin{tabular}{ |p{1.4cm}||p{1.4cm}|p{1.4cm}|p{1.4cm}|}
\hline
&$\lfloor\gamma_{CSSF}\rfloor$&$\lfloor\gamma_{HM1}\rfloor$&$\lfloor\gamma_{HM2}\rfloor$\\
\hline\hline
$\lfloor\gamma_{CSSF}\rfloor$&
& 0.0\%  & 5.0\%  \\
\hline
$\lfloor\gamma_{HM1}\rfloor$&
2.5\%   &   & 7.5\%  \\
\hline
$\lfloor\gamma_{HM2}\rfloor$&
0.0\%   & 0.0\%  &\\
\hline
\end{tabular}
\caption{Comparison of the bounds on the domination number on $\mathds{G}(n,p)$.}\label{tabGamma}
$~$\\
\end{table}

\begin{table}[!ht]\centering
\begin{tabular}{ |p{1.4cm}||p{1.4cm}|p{1.4cm}|p{1.4cm}|}
\hline
&$\lfloor\gamma_{CSSF}\rfloor$&$\lfloor\gamma_{HM1}\rfloor$&$\lfloor\gamma_{HM2}\rfloor$\\
\hline\hline
$\lfloor\gamma_{CSSF}\rfloor$&
& 0.0\%  & 5.0\%   \\
\hline
$\lfloor\gamma_{HM1}\rfloor$&
74.0\%  &  & 75.4\%  \\
\hline
$\lfloor\gamma_{HM2}\rfloor$&
0.1\%   & 0.0\%  &\\
\hline
\end{tabular}
\caption{Comparison of the bounds on the domination number on  $\mathds{B}(n,p_R,p_A)$.}\label{tabGammaBip}
$~$\\
\end{table}

In Table~\ref{tabAlpha} and Table~\ref{tabAlphaBip} the percentage of considered random graphs describes the cases where the ceiled integer of the bound on the independence number describing the row is strictly greater than the ceiled integer of the bound describing the column.

\begin{table}[!ht]\centering
\begin{tabular}{ |p{1.4cm}||p{1.4cm}|p{1.4cm}|p{1.4cm}|p{1.4cm}|}
\hline
&$\lceil\alpha_{ACL}\rceil$&$\lceil\alpha_{HR}\rceil$&$\lceil\alpha_{S}\rceil$&$\lceil\alpha_{HM}\rceil$\\
\hline\hline
$\lceil\alpha_{ACL}\rceil$&
&0.0\%& 42.2\%&75.3\%\\
\hline
$\lceil\alpha_{HR}\rceil$&
38.5\%&&78.7\%&79.8\% \\
\hline
$\lceil\alpha_{S}\rceil$&
2.2\%& <0.1\% &&51.5\%\\
\hline
$\lceil\alpha_{HM}\rceil$&
0.7\% & 0.0\%  &16.9\%&\\
\hline
\end{tabular}
\caption{Comparison of the bounds on the independence number on $\mathds{G}(n,p)$.}\label{tabAlpha}
$~$\\
\end{table}

\begin{table}[!ht]\centering
\begin{tabular}{ |p{1.4cm}||p{1.4cm}|p{1.4cm}|p{1.4cm}|p{1.4cm}|}
\hline
&$\lceil\alpha_{ACL}\rceil$&$\lceil\alpha_{HR}\rceil$&$\lceil\alpha_{S}\rceil$&$\lceil\alpha_{HM}\rceil$\\
\hline\hline
$\lceil\alpha_{ACL}\rceil$&
&69.1\%& 72.8\%&85.0\%\\
\hline
$\lceil\alpha_{HR}\rceil$&
0.8\%&&29.2\%&61.9\% \\
\hline
$\lceil\alpha_{S}\rceil$&
5.7\%& 31.3\% &&50.5\%\\
\hline
$\lceil\alpha_{HM}\rceil$&
1.0\% & 7.4\%  &18.2\%&\\
\hline
\end{tabular}
\caption{Comparison of the bounds on the independence number on $\mathds{B}(n,p_R,p_A)$.}\label{tabAlphaBip}
\end{table}

All four tables shows that our new bound compare favorably to recent ones for randomly considered graphs and do not exclusively perform well on selected particular graphs. 
The impression, which we got from the proofs, that graphs ``closer'' to bipartite graphs provide a more convenient structure for our estimations in the bounds, is approved by the experiments. Tables~\ref{tabGammaBip} and~\ref{tabAlphaBip} considering the probability spaces $\mathds{B}(n,p_R,p_A)$ show significantly stronger results of all newly presented bounds in this paper compared to calculations in $\mathds{G}(n,p)$. 

%\newpage

\end{document}